 \documentclass[11pt,a4paper]{article}
\usepackage{amsmath}
\usepackage{amsfonts}
 \usepackage{amsthm}
\theoremstyle{plain}
\newtheorem{thm}{Theorem}[section]

\newtheorem{prop}[thm]{Proposition}
\newtheorem{wn}[thm]{Corollary}
\theoremstyle{definition}
\newtheorem{rem}[thm]{Remark}

\newcommand{\varcal}[1]{\mathcal{#1}}
\newcommand{\R}{\mathbb R}

\newcommand{\Z}{\mathbb Z}
\newcommand{\Rd}{{\mathbb R^d}}

\newcommand{\SprimeRd}{%
 {\ifmmode {\varcal S}'(\Rd)%
 \else ${\varcal S}'(\Rd)$\fi}}
\newcommand{\SRd}{{%
 \ifmmode \varcal S(\Rd)%
 \else $\varcal S(\Rd)$\fi}}
\newcommand{\Sprime}{{%
 \ifmmode {\varcal S}'%
 \else ${\varcal S}'$\fi}}
\newcommand{\Ss}{{%
 \ifmmode \varcal S%
 \else $\varcal S$\fi}}

\newcommand{\T}{\mathcal T}

\newcommand{\<}{\left<}

\renewcommand{\>}{\right>}

\newcommand{\abs}[1]{\left\vert#1\right\vert} 

\newcommand{\ind}{1\mkern-7mu1}
\newcommand{\UNO}{1\mkern-7mu1}
\newcommand{\convf}{\underset{f}{\Rightarrow}}
\newcommand{\Var}{\operatorname{Var}}
\newcommand{\bPsi}{{\bf \Psi}}
\newcommand{\comment}[1]{}
\numberwithin{equation}{section}

 \title{\normalsize{\textbf{PARTICLE PICTURE INTERPRETATION OF SOME GAUSSIAN PROCESSES RELATED TO FRACTIONAL BROWNIAN MOTION
}}}
 
  \author{\small{TOMASZ BOJDECKI} $^1$\\
  \normalsize{ Institute of Mathematics,
  University of Warsaw}\\
  \normalsize{ul. Banacha 2,
  02-097 Warszawa,
  Poland}\\
  \normalsize{
  e-mail: tobojd@mimuw.edu.pl}
   \and
 \small{ANNA TALARCZYK}$^{1,*}$\\
  \normalsize{ Institute of Mathematics,
  University of Warsaw}\\
  \normalsize{ul. Banacha 2,
  02-097 Warszawa,
  Poland}\\
  {\normalsize{
  e-mail: annatal@mimuw.edu.pl}}
}
\date{}
\footnotetext{$\!\!\!\!\!\!\!\!\!\!^{1}$ 
Supported in part by MNiSzW grant N N201 397537 (Poland).\\
$\!\!\!\!\!\!\!\!\!\!^{*}$ Corresponding author
}

\begin{document}
 \maketitle
  \begin{abstract}
We construct fractional Brownian motion (fBm), sub-fractional Brownian motion (sub-fBm), negative sub-fractional Brownian motion (nsfBm) and the odd part of fBm in the sense of Dzhaparidze and van Zanten (2004) by means of limiting procedures applied to some particle systems. These processes are obtained for full ranges of Hurst parameter. Particle picture interpretations of sub-fBm and nsfBm were known earlier (using a different approach) for narrow ranges of parameters; the odd part of fBm process had not been given any physical interpretation at all. 

Our approach consists in representing these processes as $\<X(1),\ind_{[0,t]}\>$, $\<X(1),\ind_{[0,t]}-\ind_{[-t,0]}\>$, $\<X(1),\ind_{[-t,t]}\>$, respectively, where $X(1)$ is an (extended) $\Sprime$-random variable obtained as the fluctuation limit of either empirical process or the occupation time process of an appropriate particle system. In fact, our construction is more general, permitting to obtain some new Gaussian processes, as well as multidimensional random fields. In particular, we generalize and presumably simplify some results by Hambly and Jones (2007). We also obtain a new class of $\Sprime$-valued density processes, containing as a particular case the density process of Martin-L\"of (1976).
\\
\vglue .3cm
\noindent

 \end{abstract}
 \newpage
 
 \section{Introduction}
  \label{sec:i}
The main objective of the present paper is to give a particle picture interpretation of the continuous
centered Gaussian processes $\xi^H=(\xi^H_t)_{t\ge 0}$,  $\zeta^H=(\zeta^H_t)_{t\ge 0}$ and $\eta^H=(\eta^H_t)_{t\ge 0}$
with  covariance functions, for $0\le s\le t$,
\begin{align}
E\xi^H_s\xi^H_t=&\frac 12 (s^{2H}+t^{2H}-\abs{t-s}^{2H}), \qquad H\in (0,1)
\label{e:1.1}\\
E\zeta^H_s\zeta^H_t=&(1-H) \left(s^{2H}+t^{2H}-\frac 12((s+t)^{2H}+\abs{t-s}^{2H})\right),
\notag\\
&\hskip 4.5cm  \ H\in (0,2)\label{e:1.2} \\
E\eta^H_s\eta^H_t=&(s+t)^{2H}-\abs{t-s}^{2H}, \qquad H\in (0,1).
\label{e:1.3}
\end{align}
$\xi^H$ is the well-known fractional Brownian motion (fBm), $\zeta^H$, for $H<1$, is the sub-fractional Brownian motion (subfBm, see e.g. \cite{sub_fbm}), and for $H>1$ is the so-called negative sub-fractional Brownian motion (nsfBm \cite{extfbm}), $\eta^H$ is a counterpart of the subfBm in the  sense that both processes occur in the decomposition of fBm studied in \cite{DzvZ}.
 It is called the odd part of fBm in that paper. It is also related to nsfBm, namely, $\zeta_t^H=K\int_0^t\eta_s^{H-1}ds$ for a constant $K$ (see \cite{extfbm}).  It is known that the closures of the intervals written above are maximal ranges of the Hurst parameters $H$.  In each case, the extreme points of these intervals correspond to trivial processes, which will be excluded from our considerations.
 
There are many models, related to particle systems or not, which lead to fBm (see, e.g., \cite{functlim1},\cite{MR0400329}, \cite{MR2477392}, \cite{Enriquez_fBm}), especially, for $H\ge \frac 12$. Sub-fBm and nsfBm appeared in a natural way in connection with occupation time fluctuations of particle systems (\cite{Deuschel1994}, \cite{BirknerZahle}, \cite{functlim1}, \cite{functlim6}, \cite{functlim8}, we give some more details below) and this was the reason for us to study them in \cite{sub_fbm} and \cite{extfbm}. Sub-fractional Brownian motion has gained independent interest and was investigated also by other authors (e.g. \cite{Bardina_Bascompte}, \cite{Tudor_subfBm}, \cite{MR2506510}, \cite{MR2593566}). However, the particle models studied earlier led to these processes for narrow ranges of parameters only; namely, $H\in[\frac 12,1)$ for sub-fBm and $H\in (1, \frac 54]$ for nsfBm. It seemed natural, and interesting, to ask if there exist other particle models that would permit to obtain these processes for the whole ranges of parameters. For example, in the case of sub-fBm, the difficulty of interpretation seemed to follow from the fact that for $H<\frac 12$ increments of this process on non-overlapping intervals are negatively correlated. As far as we know, the process $\eta^H$ has not been given any physical interpretation at all.

In this paper we use ``the white-noise approach'', i.e., our starting point is the well-known fact that the standard Brownian motion can be represented as $\<X,\ind_{[0,t]}\>$, where $X$ is the white noise (see, e.g. \cite{Kuo}). A similar construction permits to obtain processes defined in $\eqref{e:1.1}-\eqref{e:1.3}$ for all possible $H$. It also yields new interpretations of sub-fBm and nsfBm for parameters obtained before. 

The second objective, related to the first one, is to introduce a new class of density processes in $\Sprime$, and associate it to a class of real Gaussian processes, containing as a special case the process studied by Hambly and Jones (\cite{hambly_jones}, \cite{hambly_jones_errata}).

All the models discussed in this paper are based on the following particle system in $\R$: at time $t=0$ positions of the particles are determined by a point measure $\nu$ (in general random), which is in some sense homogeneous; then they are evolving independently according to the standard, symmetric $\alpha$-stable L\'evy motion ($\alpha \in (0,2]$). In some models they additionally undergo a critical binary branching. The evolution of the system is described by the empirical process $N=(N_t)_{t\ge 0}$, where $N_t(A)$ is the number of particles in the set $A\subset\R$ at time $t$. In general, the system depends on a parameter $T$ which will tend to infinity. The corresponding empirical process will be denoted by $N^T$.

We consider three classes of models. The first one gives a particle picture interpretation of $\xi^H$, $\zeta^H$ and $\eta^H$ for $H<\frac 12$ in a somewhat circuitous way (two limit passages), but, on the other hand, it is quite general and leads to a new class of the so called density processes. Moreover, we obtain a convergence result of Hambly and Jones \cite{hambly_jones}, \cite{hambly_jones_errata}, as a very 
particular case.
The initial configuration is given by a measure $\nu_T$ such that as $T\to \infty$ the density of particles increases. 

We define
\begin{equation}
 X_T(t)=\frac{N_t^T-EN_t^T}{\sqrt{T}}.
 \label{e:1.4}
\end{equation}
If $\nu_T$ is a Poisson random measure with intensity measure $T\lambda$, where $\lambda$ is the Lebesgue measure, it is well known that $X_T$, regarded as $\Sprime$-valued processes ($\Sprime$ is the space of tempered distributions) converge in law to the so-called density process, which is a continuous centered Gaussian $\Sprime$-valued process of Ornstein-Uhlenbeck type (see e.g. \cite{MartinLof},\cite{MR527199},\cite{MR690135}).  Our first result is a generalization of this fact (in the non-branching case). We consider a class of initial configurations $\nu_T$ which includes as particular cases both homogeneous Poisson random measure and the deterministic measure $\sum_{j\in \Z}\delta_{\frac jT}$. We show convergence of finite dimensional distributions in $\Sprime$ to a certain new density process $X$ (Theorem \ref{thm:2.1}). Moreover, we show that this convergence holds for a class of test functions wider than $\Ss$ (the space of smooth rapidly decreasing functions).
 We then study real processes of the form 
 \begin{equation}
  (\<X(1),\psi_t\>)_{t\ge 0},
  \label{e:1.5}
 \end{equation}
where $\<X(1),\varphi\>$ is defined by a natural extension.

The most important examples are 
\begin{align}
 \psi_t=&\ind_{[0,t]}\label{e:1.6}\\
 \psi_t=&\ind_{[0,t]}-\ind_{[-t,0]}\label{e:1.7}\\
 \psi_t=&\ind_{[-t,t]}\label{e:1.8}
\end{align}
Note that the space parameter plays now the role of ``time''. For $\psi_t$ of the form \eqref{e:1.6} we obtain a generalization of a convergence result by Hambly and Jones (\cite{hambly_jones}, \cite{hambly_jones_errata}) who took a  diffferent approach and did not use explicitely the high density of the system (see Remark \ref{rem:2.4} (b)). Note that the processes $\<N_1^T,\psi_t\>$ in each case $\eqref{e:1.6}-\eqref{e:1.8}$ have a clear physical interpretation. For \eqref{e:1.7} the interpretation of \eqref{e:1.4} is particularly nice, since $E\<N_1^T,\psi_t\>=0$ if $\nu_T$ is symmetric.

Next, assuming that the initial configuration is deterministic or nearly deterministic, we make a simple second passage to the limit. We speed up the ``time'' in \eqref{e:1.5} (in fact, shrinking the space), i.e. we consider $\<X(1),\psi_{Tt}\>$, 
or more generally, $\<X(1),\varphi(\frac \cdot T)\>$
and normalize appropriately. 
In the limit, as  $T\to \infty$, for $\psi_t$ of the form \eqref{e:1.6} we get easily a fractional Brownian motion (this result was proved by Hambly and Jones \cite{hambly_jones}), for \eqref{e:1.7} a sub-fractional Brownian motion, and for \eqref{e:1.8} the process $\eta^H$ defined in \eqref{e:1.3}, with $H<\frac 12$ in each case (Propositions \ref{prop:2.5}, \ref{prop:2.6}). In this way we obtain the desired particle picture interpretation of the sub-fBm for small Hurst parameters. We remark that this procedure for truly random initial measures does not lead to interesting results.

It is worthwhile to note that writing the processes in the form \eqref{e:1.5} permits to derive easily their properties, such as stationarity of increments, negative correlatedeness of increments, long range dependence. Some of these properties are preserved after the second passage to the limit.
 
In the second  model, giving the particle picture interpretation of $\xi^H$, $\zeta^H$, $\eta^H$ for $H<\frac 12$ more directly, we consider an initial configuration determined by a fixed measure $\nu$, which is deterministic, or almost deterministic. We assume $\alpha<1$; again, there is no branching. We carry out the same space transformation as in the second passage to the limit in the previous model and define
\begin{equation}
\<Z_T,\varphi\>=\frac{\<N_1,\varphi(\frac \cdot {T})\>-E\<N_1,\varphi(\frac \cdot T)\>}{\sqrt{T^{1-\alpha}}}.
 \label{e:1.9a}
\end{equation}
We prove convergence in law in $\Sprime$ to a centered Gaussian $\Sprime$-variable. Again, this result can be extended to test functions $\psi_t$  of the forms \eqref{e:1.6}-\eqref{e:1.8}, yielding convergence of $\<Z_T,\psi_\cdot\>$ in the sense of finite dimensional distributions to fBm, sub-fBm, $\eta^H$, respectively, with $H<\frac 12$ (Theorem \ref{thm:2.9}). 

The third class of models is related to occupation time fluctuations. We consider a particle system with homogeneous Poisson initial condition with or without high density, with or without branching and define 
\begin{equation}
 Y_T(t)=\frac 1{F_T}\left(\int_0^{Tt}N_s^Tds-\int_0^{Tt}EN_s^Tds\right).
 \label{e:1.9}
\end{equation}
where $F_T$ is a deterministic norming.

In \cite{functlim1} and \cite{functlim6}, for ``large'' $\alpha$ ($\alpha>1$ in the case without branching, and $\alpha>\frac 12$ in the case with branching), we showed the weak functional convergence of $Y_T$ in $\Sprime$, as  $T\to \infty$,  to a limit process $Y$ which
is of the form $K\vartheta \lambda$, where  $\vartheta$ is a centered Gaussian process. Depending on the specific parameters of the model, $\vartheta$ is 
\begin{itemize}
 \item a fractional Brownian motion with $H\in(\frac 12,\frac 34] $,
 \item a subfractional Brownian motion with $H\in(\frac 12,1)$,
 \item a negative subfractional Brownian motion with $H\in (1,\frac 54)$.
\end{itemize}
Considering a branching system in equilibrium Mi\l o\'s\ in \cite{Milos1} obtained a fractional Brownian motion for $H\in (\frac 12,1)$. On the other hand, for ``small'' $\alpha$ the limit $Y$ also exists (\cite{functlim2}) but is of a completely different form. It is truly $\Sprime$-valued process with independent increments. 

In this paper, in the non-branching case with $\alpha<1$,  we  consider
$\<Y_T(1),\psi_t\>$ for $\psi_t$ of the form \eqref{e:1.6}-\eqref{e:1.8}, similarly as for density processes,  and in the limit we obtain fBm, sub-fBm and the process
$\eta^H$ with $H\in(\frac 12,1)$. The fact that the formal expression $\<Y(1),\ind_{[0,t]}\>$ leads to a fBm was observed in \cite{1104_0056v1}. Moreover, for $\alpha>1$, considering $\<Y_T,\psi_t\>$ with $\psi_t$ of the form \eqref{e:1.7} we get a nsfBm with $H\in (1,\frac 32)$. In this case the
 norming $F_T$ is different from the one in \cite{functlim1}, due to the fact that $\int_\R\psi_t(x)dx=0$.
 
 So we see that as $\alpha$ changes  from ``small'' to ``large'' values, there occurrs a sort of phase transition in the sense that a complicated time structure and simple space structure changes into a simple temporal structure and complicated spatial structure. Moreover, this ``complicated'' structure in both cases is related to fractional Brownian motion.
 
 It remains to find a particle picture yielding a negative sub-fractional Brownian motion for the missing range of parameters $H\in[\frac 32,2)$. This is achieved by considering a high density branching particle system with high branching intensity
(Theorem \ref{thm:2.15}).

We remark that the argument of this paper can be carried out in the multidimensional case, for example we can consider $\psi_{t_1,\ldots, t_d}=\ind_{[0,t_1]\times\ldots\times [0,t_d]}$, obtaining some Gaussian random fields. 

It seems more difficult to try to extend the approach of this paper to the systems with infinite variance branching. It is worthwhile to notice that taking formally $\<X(1),\ind_{[0,t]}\>$, where $X$ is the stable $\Sprime$-valued limit process in Theorem 2.1 in \cite{functlim4}, one obtains a linear fractional stable motion with parameters $a=1, b=-1, \frac 1{2}< H<1$ defined in \cite{Samor_Taqqu} (Definition 7.4.1). One would obtain an interesting and probably new interpretation of this process if this procedure were justified by a suitable limit theorem.

The following notation is used in the paper:\\
$\lambda$: Lebesgue measure;\\
$\Ss$: space of $C^\infty$ rapidly decreasing functions on $\R$;\\
$\Sprime$: space of tempered distributions (topological dual of $\Ss$);\\
$\<\cdot, \cdot\>$: duality, in particular on $\Sprime\times \Ss$ and $\<\mu,f\>=\int fd\mu$;\\
$\convf$: weak convergence of finite dimensional distributions of processes in appropriate space;\\
$p_t(x)$: transition probability density of the standard symmetric $\alpha$-stable L\'evy process in $\R$;\\
$\T_t$: semigroup determined by $p_t$, i.e. $\T_t\varphi=p_t*\varphi$;\\
$\hat \varphi(x)=\int_\R e^{ixy}\varphi (y)dy$.\\
Generic constants are denoted by $C, C_i$ with possible dependencies in parenthesis or in subscripts.

In Section 2 we describe the particle systems, formulate the results and discuss them. Section $3$ contains proofs.

\section{Results}
\label{sec:results}

\subsection{Density processes}
\label{sec:density}

We consider the following particle system introduced in \cite{functlim8}.

Let $\theta$ be a non-negative integer-valued random variable with distribution
\begin{equation}
\label{e:2.1}
P(\theta=k)=p_k,\,\, k=0,1,2,\ldots,
\end{equation}
and such that $E\theta^3<\infty$. 
Let $\theta_j,j\in\Z$, be independent copies of $\theta$, and for each $j\in\Z$ and $k=1,2,\ldots$, let $\rho^j_k=(\rho^j_{k,1},\ldots,\rho^j_{k,k})$ be a random vector with values in $[j,j+1)^k$.
We assume that 
$(\theta_j,(\rho^j_k)_{k=1,2},\ldots)$, $j\in \Z$, are independent.
Given $T>0$, these  objects determine a random point measure $\nu_T$ on $\R$ in the following way: For each $j, \theta_j$ is the number of points in the interval $[\frac jT,\frac{j+1}T)$, and for each $k$, if $\theta_j=k$, the positions of those points are determined by $\frac 1T\rho^j_k$. In other words,
\begin{equation}
\label{e:2.2}
\nu_T=\sum_{j\in\Z}\sum^{\theta_j}_{n=1}
\delta_{\kappa_{j,n,T}},
\end{equation}
where 
\begin{equation}
\label{e:2.2a}
\kappa_{j,n,T}=\frac 1T\rho^j_{\theta_{j},n},
\end{equation}
 and $\delta_a$ is the Dirac  measure at $a\in\R$.

Observe that both the deterministic measure 
\begin{equation}
\nu_T=\sum_{j\in \Z} \delta_{\frac j T}
 \label{e:2.3}
\end{equation}
and the Poisson random measure with intensity $T\lambda$ are of the form \eqref{e:2.2}. 

Fix $\alpha\in(0,2]$ and assume that at the initial time $t=0$ there is a collection of particles in $\R$ with positions determined by a measure $\nu_T$ of the form \eqref{e:2.2}. As time evolves, these particles move independently according to the symmetric $\alpha$-stable L\'evy process. 

We define a signed measure-valued process $X_T$ by \eqref{e:1.4}, for which we have the following result.

\begin{thm}
 \label{thm:2.1}
 $X_T\convf X$ in $\Sprime$ as $T\to\infty$, where $X$ is a centered Gaussian $\Sprime$-valued process with covariance functional
 \begin{equation}
 E\<X(s),\varphi\>\<X(t),\psi\>=E\theta\int_\R \varphi(x)\T_{\abs{t-s}}\psi(x)dx
 +(\Var \theta -E\theta) \int_\R \T_s \varphi(x)\T_t \psi(x)dx.
  \label{e:2.4}
 \end{equation}
 \end{thm}
 
\begin{rem}
\label{rem:2.2}
(a) The limit process $X$ will be called density process.  Note that if $\nu_T$ is a Poisson random measure with intensity $T\lambda$ or, more generally, if $\theta$ is such that $E\theta=\Var \theta$, then $X$ is the classical density process (e.g. \cite{MartinLof}).\\
(b) The same result is true if the $\alpha$-stable motion is replaced by a more general process $\vartheta$, with the covariance functional of the limit written as 
 \begin{equation}
 E\theta\int_\R 
 E\varphi(x + \vartheta_s)\psi(x+\vartheta_t)dx
 +(\Var \theta -E\theta) \int_\R  E\varphi(x+\vartheta_s)E\psi(x+\vartheta_t)dx.
  \label{e:2.5}
 \end{equation}
 The only condition on $\vartheta$ is that there exists $m\ge 0$ such that \eqref{e:2.5} is finite if $\varphi$ and $\psi$ are replaced by $\phi_m(x)=\frac 1{1+\abs x^m}$.\\
 (c) There is no problem to extend this result to $\R^d$. In the definition of $\nu_T$ the interval $[j,j+1)$ should be replaced by the cube $[j_1,j_1+1)\times\ldots\times[j_d,j_d+1)$ and in \eqref{e:1.4} the normalization is $\sqrt{T^d}$. The limit is an $\SprimeRd$-valued process with covariance of the form \eqref{e:2.5}.
\end{rem}

We want to extend $\<X,\psi\>$ for some $\psi$'s  which are not necessarily in $\Ss$. Let $\bPsi$ be the class of bounded piecewise continuous functions with compact support. From \eqref{e:2.4} it is clear that for $\psi\in\bPsi$ we can define $\<X,\psi\>$ by 
$L^2$-approximation. Let $\Psi=(\psi_t)_{t\ge 0}$ be a family of functions from $\bPsi$. Typically, $\psi_t$ are of the form \eqref{e:1.6}-\eqref{e:1.8}. We will study (centered Gaussian) processes of the form 
\begin{equation}
\varrho_t^\Psi=\<X(1),\psi_t\>_{t\ge 0}.
 \label{e:2.6}
\end{equation}
By \eqref{e:2.4}, the covariance function of $\varrho^\Psi$ is 
\begin{equation}
E\varrho_s^\Psi\varrho_t^\Psi=E\theta \int_\R \psi_s(x)\psi_t(x)dx
 +(\Var \theta -E\theta) \int_\R \psi_s(x)\T_2 \psi_t(x)dx.
 \label{e:2.7}
\end{equation}

The next proposition shows that this extension is compatible with Theorem \ref{thm:2.1}.

\begin{prop}
 \label{prop:2.3}
 \begin{equation*}
 \left(\<X_T(1),\psi_t\>\right)_{t\ge 0}\convf \varrho^\Psi, \qquad \textrm{as}\ T\to \infty.
 \end{equation*}
\end{prop}

\begin{rem}
\label{rem:2.4}
(a) This proposition gives a particle picture interpretation of the process $\varrho^\Psi$.
\medskip
 
\noindent
(b) For deterministic $\theta$ and $\Psi$ of the form \eqref{e:1.6}, the process $\varrho^\Psi$ is, up to a constant, the same as the process $G$ defined in Proposition 4.2 of \cite{hambly_jones}, (see also \cite{hambly_jones_errata}) for parameter $c=1$. Moreover, the convergence result of that proposition follows directly from Proposition \ref{prop:2.3} by a standard conditioning. (The   extension to general $c$ is also immediate.) We stress that the proof in \cite{hambly_jones} is completely different, not exhibiting the high density of the system.
\end{rem}

The covariance \eqref{e:2.7} may look complicated but it does permit to deduce some properties of the underlying process (see Proposition \ref{prop:2.7}). The situation simplifies considerably after the second passage to the limit which consists in shrinking the space. 

For a function $f:\R\mapsto \R$ and $T>0$ we denote
\begin{equation}
 \label{e:2.8}
 S_T f(x)=f(\frac x T).
\end{equation}
We define an $\Sprime$-valued random variable $\tilde Z_T$ by
\begin{equation}
 \label{e:2.8a}
 \<\tilde Z_T, \varphi\>=\frac 1{F_T} \<X(1), S_T\varphi\>, \qquad \varphi \in \Ss,
\end{equation}
where $X$ is given by Theorem \ref{thm:2.1} and $F_T$ is a suitable norming. We have the following simple fact.

\begin{prop}
 \label{prop:2.5}
 (a) If $\Var \theta=0$ and $F_T=\sqrt{T^{1-\alpha}}$ then $\tilde Z_T$ converges in law, as $T\to \infty$, to a centered Gaussian $\Sprime$-valued random variable $Z$ with covariance 
 \begin{equation}
  \label{e:2.9}
  E\<Z,\varphi\>\<Z,\psi\>=\frac {E\theta}\pi \int_\R \hat\varphi (x)\overline{\hat \psi(x)} \abs{x}^\alpha dx. 
 \end{equation}
(b) If $\Var \theta>0$ and $F_T=\sqrt T$ then $\tilde Z_T$ converges in law in $\Sprime$ to $(\Var \theta)\, W$, where $W$ is the white noise. 
\end{prop}

From the proof it is clear that in this proposition the functions $\varphi\in \Ss$ can be replaced by $\psi_t\in \bPsi$, giving the convergence of finite dimensional distributions of the processes $\frac{1}{F_T}\varrho^{S_T\Psi}_\cdot:=\frac{1}{F_T}\<X(1),S_T\psi_\cdot\>$ (cf. \eqref{e:2.6}) provided that $\int_\R\abs{\hat\psi_t(x)}^2\abs{x}^\alpha dx<\infty$ in case (a).
In particular, if $\Psi$ has one of the forms \eqref{e:1.6}-\eqref{e:1.8}, then
\begin{equation}
 \label{e:2.10}
\varrho_t^{S_T\Psi}=\varrho_{Tt}^\Psi
\end{equation}
and we have the following Proposition:

\begin{prop}
 \label{prop:2.6}
 Assume $\Var \theta=0$ and $\alpha<1$. Let $\Psi$ have one of the forms \eqref{e:1.6}-\eqref{e:1.8}. Then the processes $(\frac 1 {\sqrt {T^{1-\alpha}}}\varrho^\Psi_{Tt})_{t\ge 0}$ converge in the sense of finite dimensional distributions to:\\
 (a) a fractional Brownian motion $K\xi^H$ if $\psi_t=\ind_{[0,t]}$,\\
 (b) a sub-fractional Brownian motion $K\zeta^H$ if $\psi_t=\ind_{[0,t]}-\ind_{[-t,0]}$,\\
 (c) the process $K\eta^H$ if $\psi_t=\ind_{[-t,t]}$,\\
with $H=\frac{1-\alpha}{2}$ in each case. 
\end{prop}

We have assumed that $\Var \theta=0$, i.e., the initial configuration of particles is close to deterministic, since for random $\theta$, in all cases \eqref{e:1.6}-\eqref{e:1.8} one obtains the Brownian motion. This proposition gives the desired particle picture interpretation for sub-fBm with $H<\frac 12$ as well as an interpretation of the process $\eta^H$. Part $(a)$ was observed in \cite{hambly_jones}. 

We close this  subsection with a brief discussion of the properties of the process $\varrho^\Psi$  and of the $\Sprime$-random variable $X(1)$.

\begin{prop}
 \label{prop:2.7}
Fix $\Psi=(\psi_t)_{t\ge 0}$, $\psi_t\in \bPsi$, and let $\varrho^\Psi$ be the process defined by \eqref{e:2.6}, and $X(1)$ be given by Theorem \ref{thm:2.1}.
\medskip
 
 \noindent
 (a) Homogeneity:
 $\<X(1),\varphi\>$ has the same distribution as $\<X(1),\varphi(\cdot -a)\>$ for $a\in\R$, $\varphi\in \Ss$, and, for $\Psi$ of the form \eqref{e:1.6}, $\varrho^\Psi$ has stationary increments.
 \medskip
 
 \noindent
 (b) Correlation of increments: Assume that if $s<t\le u<v$ then $\psi_t-\psi_s$ and $\psi_v-\psi_u$ have disjoint supports; moreover, the function $t\mapsto \psi_t(x)$ is either non-decreasing for $x\in \R$ or this property holds on $\R_+$ and the functions $\psi_t$ are odd. Then\\
 (i) if $\Var \theta<E \theta$ then the increments of $\varrho^\Psi$ are negatively correlated;\\
 (ii) if $\Var \theta=E\theta$ then $\varrho^\Psi$ has independent increments;\\
 (iii) if $\Var \theta > E \theta$ then the increments of $\varrho^\Psi$ are positively correlated.
 \medskip
 
 \noindent
 (c) Long range dependence: 
If $\varphi_1, \varphi_2$ have compact supports then
\begin{equation}
 \label{e:2.12}
 \lim_{\tau\to\infty}{\tau^{1+\alpha}}E\<X(1),\varphi_1\>\<X(1),\varphi_2(\cdot - \tau)\>
=C(\Var\theta -E\theta)\int_\R\varphi_1(x)dx\int_\R \varphi_2(x)dx.
\end{equation}
In particular, for $\Psi$ of the form \eqref{e:1.6}
 \begin{equation}
  \label{e:2.11}
  \lim_{\tau\to \infty} \tau^{1+\alpha} E\left(\varrho_t^\Psi -\varrho_s^\Psi \right) \left(\varrho_{v+\tau}^\Psi -\varrho_{u+\tau}^\Psi \right)
  =C(\Var \theta-E\theta)(t-s)(v-u).
 \end{equation}
(d) Path continuity: If 
\begin{equation}
 \label{e:1.13}
\int_\R(\psi_{t_1}(x)-\psi_{t_2}(x))^2dx\le C(T) \abs{t_1-t_2}^\beta, \quad t_1,t_2\le T,
\end{equation}
for some $2>\beta>0$ and any $T>0$,  then $\varrho^\Psi$ has a continuous version; more precisely it is locally H\"older continuous with exponent $<\frac \beta 2$.
\end{prop}

All these properties follow easily from \eqref{e:2.4}, \eqref{e:2.7}. We will give a brief explanation in the next section.

\begin{rem}
 \label{rem:2.8} (a) 
 A similar calculation as in deriving \eqref{e:2.11} shows that the long range dependence rate of increments of $\varrho^\Psi$ with $\Psi$ of the form \eqref{e:1.8} is also $\tau^{-(1+\alpha)}$. On the other hand, this rate for $\Psi$ of the form \eqref{e:1.7} is $\tau^{-(2+\alpha)}$ due to the fact that in this case $\psi_t$ are odd.\\
 (b) Properties (a) and (b) are clearly preserved after the second passage to the limit. In particular, we thus obtain a simple proof of negative correlatedness of increments of sub-fractional Brownian motion for $H<\frac 12$. This fact was proved in \cite{sub_fbm} but the argument was rather cumbersome.\\
 (c) The relationship between the sign of $\Var \theta -E \theta$ and positive/negative correlatedness of increments seems quite interesting and unexpected. In particular, we see that independence of incrementss occurs not only in the Poisson case.\\
 (d) In the special case $\Var \theta=0$ and $\Psi$ of the form \eqref{e:1.6} property \eqref{e:2.11} and stationarity of increments of $\varrho^\Psi$ was obtained in \cite{hambly_jones} (see Remark \ref{rem:2.4}(b)).\\
 (e) In general, the process $\varrho^\Psi$ is not self-similar. Nevertheless, if $F_T$ is regularly varying at infinity, then the $\convf$ limit of $\left(\frac 1{F_T}\varrho^{\Psi}_{Tt}\right)_{t\ge 0}$, if it exists, is a self-similar process (cf. \eqref{e:2.10} and Proposition \ref{prop:2.6}).\\
 (f) This proposition can be reformulated for more general particle motions and for the multidimensional case.
\end{rem}

\subsection{Particle picture for small Hurst parameters; direct approach.}
\label{sec:2.2}
As announced in the Introduction, in this subsection we show how to obtain fBm, sub-fBm, $\eta^H$ with $H<\frac 12$, employing just one passage to the limit.

Consider a particle system described in the previous subsection with initial configuration determined by a measure $\nu$ given by \eqref{e:2.2} with $T=1$ and deterministic $\theta$. Recall that this means that the initial number of particles in each interval $[j,j+1)$, $j\in \Z$, is fixed, non-random (see \eqref{e:2.1}).

\begin{thm}
 \label{thm:2.9}Assume that $\alpha<1$ and define $Z_T$ by \eqref{e:1.9a}.\\
 (a) $Z_T$ converges in law, as $T\to \infty$, in $\Sprime$ to the random variable $Z$ as in Proposition \ref{prop:2.5}(a). \\
 (b) Let $(\psi_t)_{t\ge 0}$ have one of the forms \eqref{e:1.6}-\eqref{e:1.8}. Then
 \begin{equation*}
  \left(\<Z_T,\psi_t\>\right)_{t\ge 0}\convf K\vartheta^H, \qquad \textrm{as}\ \ T\to\infty,
 \end{equation*}
where
\begin{equation}
 \label{e:2.14}
 \vartheta^H=\begin{cases}
              \xi^H \qquad &\textrm{if}\ \ \psi_t=\ind_{[0,t]},\\
               \zeta^H &\textrm{if}\ \ \psi_t=\ind_{[0,t]}-\ind_{[-t,0]},\\
                \eta^H  &\textrm{if}\ \ \psi_t=\ind_{[-t,t]},
             \end{cases}
\end{equation}
with $H=\frac{1-\alpha}{2}$.
\end{thm}
\begin{rem}
 \label{rem:2.10} Theorem \ref{thm:2.9} has an immediate extension to the $d$-dimensional case. The condition $\alpha<1$ is then replaced by $\alpha<d$.
\end{rem}

\subsection{Processes resulting from occupation limits}
\label{sec:occupation}
So far we have obtained the particle picture for the processes $\xi^H$, $\zeta^H$ and $\eta^H$ with $H<\frac 12$. In this subsection we show how our ``white noise'' approach permits to obtain them for $H>\frac 12$. 

We consider the same particle system as before, with or without branching, we assume that the initial configuration is homegeneous Poisson. We define $Y_T$ by \eqref{e:1.9} with an appropriate $F_T$. (If the system does not depend on $T$ we write $N$ instead of $N^T$).

\begin{prop}
 \label{prop:2.11}
 Assume that there is no branching, the initial measure is Poisson with intensity $\lambda$ and $\alpha<1$. Let $F_T=\sqrt{T}$ and consider $\Psi=(\psi_t)_{t\ge 0}$, $\psi_t\in\bPsi$. Then
 \begin{equation*}
  \left(\<Y_T(1),\psi_t\>\right)_{t\ge 0}\convf \vartheta \qquad \textrm{as} \ \ T\to \infty,
 \end{equation*}
where $\vartheta$ is a centered Gaussian process with covariance
\begin{equation}
 \label{e:2.15}
 E\vartheta_t\vartheta_s=\frac 1\pi \int_\R \hat\psi_t(x)\overline{\hat \psi_s(x)}\frac 1{\abs x^\alpha} dx.
\end{equation}
\end{prop}

\begin{wn}
 \label{wn:2.12}
 If $\Psi$ has one of the forms \eqref{e:1.6}-\eqref{e:1.8} then the limit $\vartheta$ has the form $\vartheta_t=K\vartheta^H_t$ where $\vartheta^H$ is as in \eqref{e:2.14} with $H=\frac{1+\alpha}{2}$.
\end{wn}

\begin{rem}
 \label{rem:2.13} (a) The covariance \eqref{e:2.15} can be written in a more explicit form 
 \begin{equation*}
  E\vartheta_t\vartheta_s=\frac{2^{1-\alpha}\Gamma(\frac{1-\alpha}{2})}{\sqrt \pi\Gamma(\frac \alpha 2)}
  \int_{\R^2}\frac{\psi_t(x)\psi_s(x)}{\abs{x-y}^{1-\alpha}}dxdy.
 \end{equation*}
(b) The process $\vartheta$ which appears in Proposition \ref{prop:2.11} can be obtained by putting formally $\<W_0^{(\alpha)}(1),\psi_t\>$, where the ${\cal S}'$-valued Wiener process $W_0^{(\alpha)}$ is the weak functional limit of $Y_T$ (see Theorem 2.1 in  \cite{functlim2}). Proposition \ref{prop:2.11}  shows that this formality is justified. The fact that 
$\<W_0^{(\alpha)}(1),\ind_{[0,t]}\>$ is a fractional Brownian motion was observed in \cite{1104_0056v1}.
\medskip

\noindent
(c) For ``large'' $\alpha$, fBm and sub-fBm with large $H$ appeared in the temporal structure of the limit of $Y_T$ in models with or without branching and with Poisson, equilibrium or deterministic initial conditions (see \cite{functlim1}, \cite{Milos1}, \cite{functlim8}). Also nsfBm with $H\in(1,\frac 54)$ was obtained in this way from a branching Poisson system with high initial density and $\alpha>1$ (see \cite{functlim6}).
\end{rem}

In the next proposition we show that the negative sub-fractional Brownian motion, and with a wider range of parameter $H$, appears also in the spatial structure of the limit.

\begin{prop}
 \label{prop:2.14}
 Assume that there is no branching, the initial measure is Poisson with intensity $\lambda$, $\alpha>1$ and $F_T=\sqrt{T}$. For $\psi_t$ of the form \eqref{e:1.7}
 \begin{equation*}
   \left(\<Y_T(1),\psi_t\>\right)_{t\ge 0}\convf K\zeta^H \qquad \textrm{as} \ \ T\to \infty,
 \end{equation*}
with $H=\frac{1+\alpha}{2}\in(1,\frac{3}{2}]$.
\end{prop}

The covariance of $\zeta^H$ has again the form \eqref{e:2.15}, where the corresponding integral is finite due to the special form of $\psi_t$. In fact, this integral is finite for $\alpha<3$. The restricted range of $H$ ($H\le \frac 32$) follows from our model, which requires $\alpha\le 2$. To obtain the particle picture interpretation for nsfBm with the full range of parameter $H$ we have to use a more complex model.

\begin{thm}
 \label{thm:2.15} Assume that $\frac 12<\alpha<\frac 32$. Consider binary branching particle system with branching rate $V_T\to \infty$. Assume that the initial measure is Poisson with intensity $H_T\lambda$ such that 
 \begin{equation}
  \label{e:2.16}
  \lim_{T\to\infty} T^7 V_T H_T^{-1}=0.
 \end{equation}
Let $F_T=\sqrt{TV_T H_T}$. Then for $\psi_t$ of the form \eqref{e:1.7},
 \begin{equation*}
   \left(\<Y_T(1),\psi_t\>\right)_{t\ge 0}\convf K\zeta^H \qquad \textrm{as} \ \ T\to \infty,
 \end{equation*}
with $H=\frac{1+2\alpha}{2}\in(1,2)$.
\end{thm}

\begin{rem}
 \label{rem:2.16} (a) The same model with $\alpha<\frac 12$ gives in the limit a sub-fBm with $\frac 12<H<1$. It seems interesting to observe that $\alpha=\frac 12$ does not give the null process (see \eqref{e:1.2}). This is due to the fact (easy to verify) that
 \begin{equation}
 \label{e:2.17}
  E\zeta_s^H\zeta_t^H=C(H)\int_\R \frac{(1-\cos(sx))(1-\cos(tx))}{\abs{x}^{1+2H}}dx
 \end{equation}
for all $H\in(0,1)\cup (1,2)$. For $H=1$, in the limit we get a non-trivial process with covariance \eqref{e:2.17}.\\
\medskip

\noindent
(b) We think that interesting results could be obtained for models studied in this subsection if branching particle systems were considered with infinite variance branching mechanism (see the remark at the end of the Introduction). This would require, however, different methods of proof.
\end{rem}

\section{Proofs}
\label{sec:proofs}

\subsection{Proof of Theorem \ref{thm:2.1}}
\label{sec:theorem2.1}

Let $\vartheta^{j,n}$, $j\in\Z$, $n=1,2,\dots$, be independent, standard $\alpha$- stable L\'evy processes ($\vartheta_0^{j,n}=0$), independent of $\nu_T$. The process  $X_T$ defined by \eqref{e:1.4} has the form

\begin{equation}
\label{e:3.1}
X_T(t)=\frac{1}{\sqrt{ T}}\sum_{j\in\Z}(N^{T,j}_t-EN^{T,j}_t),
\end{equation}
where $N^{T,j}$ is the empirical process of the system which at time $t=0$ starts from $[j/T,(j+1)/T)$, i.e.,

\begin{equation}
\label{e:3.2}
\langle N_t^{T,j},\varphi\rangle=\sum_{n=1}^{\theta_j}\varphi(\kappa_{j,n,T}+\vartheta_t^{j,n}),
\end{equation}
see \eqref{e:2.2a}. The processes $N^{T,j}$, $j\in\Z$, are independent, hence to prove the claimed convergence we can use the central limit theorem. To this end, we show first the convergence of covariances.

By independence, \eqref{e:2.1}, \eqref{e:3.1}, \eqref{e:3.2}, we have for $s\le t$, $\varphi,\psi\in{\cal S}$,

\begin{align}
E\langle&X_s^T,\varphi\rangle\langle X_t^T,\psi\rangle \notag\\
=&\frac{1}{T}\sum_{j\in\Z}\left(\sum_{k=0}^\infty p_k\, E\Bigl(\sum_{n=1}^k\sum_{m=1}^k\varphi(\kappa_{j,n,T}+\vartheta^{j,n}_s)\psi(\kappa_{j,m,T}+\vartheta^{j,m}_t)\Bigr)\right. \notag\\
&-\left.\sum_{k=0}^\infty p_k\, E\sum_{n=1}^k\varphi(\kappa_{j,n,T}+\vartheta^{j,n}_s)\, \sum_{l=0}^\infty p_l\, E\sum_{m=1}^l\psi(\kappa_{j,m,T}+\vartheta^{j,m}_t)\right).\label{e:3.3}
\end{align}

Denote $h_{k,n}(x)=\rho_{k,n}^{[x]}-x$, where $[x]$ is the largest integer $\le x$. Clearly,

\begin{equation}
\label{e:3.4}
|h_{k,n}(x)|\le 1.
\end{equation}

Splitting the first expression on the right hand side of \eqref{e:3.3} into the sum over the diagonal ($n=m$) and the rest, by independence and Markov property of $\vartheta$, we obtain

\begin{align}
E\langle X_s^T,\varphi\rangle\langle X_t^T,\psi\rangle 
=\sum_{k=0}^\infty p_k\, &\sum_{n=1}^kI_T(k,n) +\sum_{k=0}^\infty p_k\, \sum_{
\substack{n,m=1\\ n\ne m}
}^kII_T(k,n,m)\notag\\
+&\sum_{k,l=0}^\infty p_kp_l\, \sum_{n,m=1}^kIII_T(k,n;m,l),
\label{e:3.5}
\end{align}
where (substituting $x'=x/T$)

\begin{align}
\label{e:3.6}
I_T(k,n)=&\int_\R E\, {\cal T}_s(\varphi {\cal T}_{t-s}\psi)(x+\frac{h_{k,n}(Tx)}{T})dx,\\
\label{e:3.7}
II_T(k,n,m)=&\int_\R E\, {\cal T}_s\varphi(x+\frac{h_{k,n}(Tx)}{T}){\cal T}_t\psi(x+\frac{h_{k,m}(Tx)}{T})dx,\\
\label{e:3.8}
III_T(k,n;m,l)=&\int_\R E\, {\cal T}_s\varphi(x+\frac{h_{k,n}(Tx)}{T})E\, {\cal T}_t\psi(x+\frac{h_{l,m}(Tx)}{T})dx.
\end{align}

By \eqref{e:3.4}, the terms under the integrals in \eqref{e:3.6} -- \eqref{e:3.8} converge pointwise as $T\to\infty$; on the other hand, if we denote $\phi(x)=\frac{1}{1+x^2}$, then it is easy to see that

\begin{equation}
\label{e:3.9}
{\cal T}_{t_1}(\varphi_1{\cal T}_{t_2}(\varphi_2\dots {\cal T}_{t_n}\varphi_n)\dots)(x+y) \le
C\, {\cal T}_{t_1}(\phi{\cal T}_{t_2}(\phi\dots {\cal T}_{t_n}\phi)\dots)(x)
\end{equation}
for $|y|\le 1$, $x\in\R$, $n=1,2,\dots$, provided that $|\varphi_i|\le C_1\phi$. Hence 

$$
\lim_{T\to\infty}I_T(k,n)=\int_\R{\cal T}_s(\varphi{\cal T}_{t-s}\psi)(x)dx=\int_\R\varphi(x){\cal T}_{t-s}\psi(x)dx,
$$
$$
\lim_{T\to\infty}II_T(k,n,m)=\lim_{T\to\infty}III_T(k,n;m,l)=\int_\R{\cal T}_s\varphi(x){\cal T}_t\psi(x)dx.
$$
By \eqref{e:3.5}, this shows that $E\langle X_s^T,\varphi\rangle\langle X_t^T,\psi\rangle $ converges to the right hand side of \eqref{e:2.4}.

It remains to prove that $\sum_{i=1}^ma_i\langle X_T(t_i),\varphi_i\rangle$ converges in law, as $T\to\infty$, to a Gaussian variable for any $a_1,\dots,a_m\in\R$, $t_1,\dots,t_m\ge 0$, $\varphi_1,\dots,\varphi_m\in{\cal S}$, $m=1,2,\dots$. To this end, by \eqref{e:3.1}, we apply the Lyapunov criterion. It is easy to see that the proof will be completed if we show that for any $t>0$, $\varphi\in{\cal S}$, $\varphi\ge 0$,

\begin{equation}
\label{e:3.10}
\lim_{T\to\infty}A_T(t,\varphi)=0,
\end{equation}
where 
\begin{equation}
A_T(t,\varphi)= \sum_{j\in\Z}\frac{E\langle N_t^{T,j},\varphi\rangle^3}{T^{3/2}}
=\frac{1}{T^{3/2}}\sum_{j\in\Z}\sum_{k=0}^\infty p_kE\Bigl( \sum_{n=1}^k\varphi(\kappa_{j,n,T}+\vartheta_t^{j,n})\Bigr)^3\label{e:3.10a}
\end{equation}
(see \eqref{e:3.2}, \eqref{e:2.1}, \eqref{e:2.2a}).

Using the obvious inequality $(a_1+\dots+a_k)^3\le k^2(a_1^3+\dots+a_k^3)$ for $a_1,\dots,a_k\ge0$, and employing again the functions $h_{k,n}$ we obtain

\begin{equation}
A_T(t,\varphi)\le \frac{1}{T^{3/2}}\sum_{k=0}^\infty p_kk^2\sum_{n=1}^k\int_\R{\cal T}_t\varphi^3\bigl(\frac{x}{T}+\frac{h_{k,n}(x)}{T}\bigr)dx
\le \frac{C}{\sqrt T}E\theta^3\int_\R\phi(x)dx \label{e:3.10b}
\end{equation}
by \eqref{e:3.4} and \eqref{e:3.9}, hence \eqref{e:3.10} follows.
\qed

\subsection{Proof of Proposition \ref{prop:2.3}}
\label{sec:propos2.3}
Looking at the proof of Theorem \ref{thm:2.1} it is clear that nothing changes if functions $\varphi\in{\cal S}$ are replaced by $\psi_t\in\bPsi$ (smoothness of $\varphi$ is not used there).
\qed

\subsection{Proof of Propositions \ref{prop:2.5} and \ref{prop:2.6}}
\label{sec:propos2.5.6}
As ${\tilde Z}_T$ (and $\varrho^{S_T\Psi}$) are centered Gaussian, it suffices to prove convergence of the covariances. By \eqref{e:2.4} and \eqref{e:2.8a}, we have for $\varphi,\psi\in{\cal S}$

\begin{equation}
\label{e:3.11}
E\langle {\tilde Z}_T,\varphi\rangle\langle{\tilde Z}_T,\psi\rangle= (E\theta)I_T +(\Var\theta)II_T,
\end{equation}
where
$$
I_T=\frac{1}{F_T^2}\left(\,\int_\R S_T\varphi(x)S_T\psi(x)dx-\int_\R S_T\varphi(x)\T_2S_T\psi(x)dx\right),
$$
$$
II_T=\int_\R S_T\varphi(x)\T_2S_T\psi(x)dx.
$$
Observe that $\widehat{S_T\varphi}(x)=T\hat{\varphi}(Tx)$ (see \eqref{e:2.8}) and recall that $\widehat{\T_t\varphi}(x)=e^{-t|x|^\alpha}\hat{\varphi}(x)$. Hence, by Plancherel's identity, after an obvious substitution, we have
$$
I_T=\frac{T}{F_T^2}\frac{1}{2\pi}\int_\R (1-e^{-2|x|^\alpha T^{-\alpha}})\hat{\varphi}(x)\overline{\hat{\psi}(x)}dx.
$$
This and \eqref{e:3.11} imply that if $\Var\theta=0$ and $F_T=\sqrt{T^{1-\alpha}}$, then

\begin{equation}
\label{e:3.12}
\lim_{T\to\infty}E\langle\tilde{Z}_T,\varphi\rangle\langle\tilde{Z}_T,\psi\rangle=\frac{E\theta}{\pi}\int_\R\hat{\varphi}(x)\overline{\hat{\psi}(x)}|x|^\alpha dx.
\end{equation}
On the other hand, if $\Var\theta>0$ and $F_T=\sqrt{T}$, then $\lim_{T\to\infty}I_T=0$ and the co\-variance functional of $\tilde{Z}_T$ tends to $\Var\theta\frac{1}{2\pi}\int_\R\hat{\varphi}(x)\overline{\hat{\psi}(x)}dx=\Var\theta\int_\R\varphi(x)\psi(x)dx$.

This proves Proposition \ref{prop:2.5}.

Furthermore, if $\Var\theta=0$, $\alpha<1$, and $\Psi=(\psi_t)_{t\ge 0}$ have one of the forms \eqref{e:1.6} - \eqref{e:1.8}, then it is easy to see that the same argument can be repeated with $\varphi=\psi_t$, $\psi=\psi_s$ (and  $F_T=\sqrt{T^{1-\alpha}}$) and we obtain
\begin{equation}
\label{e:3.13}
\lim_{T\to\infty}\frac{1}{T^{1-\alpha}}E\varrho^\Psi_{Tt}\,\varrho^\Psi_{Ts}= \frac{E\theta}{\pi}\int_\R\hat{\psi}_t(x)\overline{\hat{\psi}_s(x)}|x|^\alpha dx
\end{equation}
by \eqref{e:3.12}. 

It is not hard to calculate the right hand side of \eqref{e:3.13} explicitely. Depending on the specific form of $(\psi_t)_{t\ge 0}$ it is, up to a multiplicative constant, the same as the right hand side of \eqref{e:1.1}, \eqref{e:1.2}, \eqref{e:1.3}, respectively. To derive this, it is convenient to use the formula
$$
|t|^{1-\alpha}=C_\alpha\int_\R\frac{1-\cos(xt)}{|x|^{2-\alpha}}dx,
$$
valid for $\alpha<1$. For example, for $\psi_t$ of the form \eqref{e:1.7} (recall that $\psi_t=\UNO_{[0,t]}-\UNO_{[-t,0]}$) we have

\begin{align}
&\int_\R\hat{\psi}_t(x)\overline{\hat{\psi}_s(x)}|x|^\alpha dx =4\int_\R\frac{(1-\cos(xs))(1-\cos(xt))}{|x|^{2-\alpha}}dx\notag\\
&=4\int_\R\frac{(1-\cos(xs))+(1-\cos(xt))-\frac{1}{2}(1-\cos(x(s-t))-\frac{1}{2}(1-\cos(x(s+t))}{|x|^{2-\alpha}}dx\notag\\
&=4C_\alpha\Bigl(s^{1-\alpha}+t^{1-\alpha}-\frac{1}{2}(|s-t|^{1-\alpha}+(s+t)^{1-\alpha})\Bigr)\notag
\end{align}
(cf.\ \eqref{e:1.2}). Note that this is a part of the proof of \eqref{e:2.17}.
\qed

\subsection{Proof of Proposition \ref{prop:2.7}}
\label{sec:propos2.7}
Denote $\varphi^a(x)=\varphi(x-a)$. (a) follows immediately from \eqref{e:2.4} since $\T_1(\varphi^a)=(\T_1\varphi)^a$, and for $\Psi$
of the form \eqref{e:1.6} $\psi_{t+a}-\psi_{s+a}=\psi_t^a-\psi_s^a$.

(b) is a consequence of \eqref{e:2.6}, \eqref{e:2.7} and of the fact that under assumptions on $\Psi$
\begin{equation*}
 \int_\R (\psi_v(x)-\psi_u(x))\T_2 (\psi_t-\psi_s)(x)dx\ge 0.
\end{equation*}
This is obvious if $t\mapsto\psi_t(x)$ is nondecreasing for each $x\in\R$, and if this property holds on $\R_+$ and $\psi_t$ is odd for $t\ge 0$, then, by symmetry of $p_2$ we can write this integral as 
\begin{equation*}
 2\int_0^\infty\int_0^\infty (p_2(x-y)-p_2(x+y))(\psi_t(y)-\psi_s(y))(\psi_v(x)-\psi_u(x))dydx,
\end{equation*}
which is non-negative since $p_2$ is unimodal.

(c) follows from \eqref{e:2.4}, \eqref{e:2.7} ($\varphi_1\varphi_2^\tau\equiv 0$ if $\varphi_1, \varphi_2$ have compact supports and $\tau$ is large) using the well known fact that $\lim_{\tau\to \infty}\tau^{1+\alpha}p_2(x)=c$ uniformly in $x$ on compact sets.

Finally, to obtain (d), it suffices to observe that 
\begin{equation*}
 E\<X(1),\varphi\>^2\le C\int_\R \varphi^2(x)dx,
\end{equation*}
since $\int_\R p_2^2(x)dx<\infty$.
\qed

\subsection{Proof of Theorem \ref{thm:2.9}}
\label{sec:propos2.9}
For brevity, we will assume that at the beginning each interval $[j,j+1), j\in \Z$, contains exactly one particle. We then modify the notation of Section \ref{sec:theorem2.1}, writing $\vartheta^j$ instead of $\vartheta^{j,1}$, $\rho^j$ instead of $\rho^j_{1,1}$, and $h(x)=\rho^{[x]}-x$ instead of $h_{1,1}$.

Again, we use the central limit theorem, showing first the convergence of covariances and then applying the Lyapunov criterion.

(a) Fix $\varphi,\psi\in \Ss$. Note that by the self-similarity property of the $\alpha$-stable motion, $\frac {\vartheta^j_1}T$ has the same law as $\vartheta_{T^{-\alpha}}^j$. Then, similarly as in \ref{sec:theorem2.1} (c.f. \eqref{e:3.3}, \eqref{e:3.5}-\eqref{e:3.8}), we can write
\begin{align}
 E&\<Z_T,\varphi\>\<Z_T,\psi\>\notag\\
 =&\frac 1{T^{1-\alpha}}\sum_{j\in\Z} E\left(\varphi(\frac{\rho^j}{T}+\frac{\vartheta_1^j}{T})
 \psi(\frac{\rho^j}{T}+\frac{\vartheta_1^j}{T})\right)
-\frac 1{T^{1-\alpha}}\sum_{j\in\Z} E\varphi(\frac{\rho^j}{T}+\frac{\vartheta_1^j}{T})
E\psi(\frac{\rho^j}T+\frac{\vartheta_1^j}{T})\notag\\
=&\frac 1{T^{1-\alpha}}E\int_{\R^2}\varphi(\frac{x+h(x)}T+y)\psi(\frac{x+h(x)}T+y)p_{T^{-\alpha}}(y)dydx\notag\\
&-\frac 1{T^{1-\alpha}}\int_{\R^3}E\varphi(\frac{x+h(x)}T+y)E\psi(\frac{x+h(x)}T+z)p_{T^{-\alpha}}(y)p_{T^{-\alpha}}(z)dy dz dx\notag\\
=&B_T+U_T-V_T,
\label{e:3.14}
\end{align}
where
\begin{align}
 B_T=&\frac 1{T^{1-\alpha}}\left(\int_{\R^2}\varphi(\frac{x}T+y)\psi(\frac{x}T+y)p_{T^{-\alpha}}(y)dydx\right.\notag\\
 &\left. -\int_{\R^3}\varphi(\frac{x}T+y)\psi(\frac{x}T+z)p_{T^{-\alpha}}(y)p_{T^{-\alpha}}(z)dy dz dx
 \right),\label{e:3.15}\\
 U_T=&\frac 1{T^{1-\alpha}}E\int_{\R^2}\left(\varphi(\frac{x+h(x)}T+y)\psi(\frac{x+h(x)}T+y)
-\varphi(\frac{x}T+y)\psi(\frac{x}T+y)
\right)\notag\\
&\hskip 7cm p_{T^{-\alpha}}(y)dydx,\label{e:3.16}\\
V_t=&\frac{1}{T^{1-\alpha}}
\int_{\R^3}\left(E\varphi(\frac{x+h(x)}T+y)E\psi(\frac{x+h(x)}T+z)-\varphi(\frac{x}T+y)\psi(\frac{x}T+z)\right)\notag\\
&\hskip 6cm p_{T^{-\alpha}}(y)p_{T^{-\alpha}}(z)dy dz dx.
\label{e:3.17}
\end{align}
In $B_T$ we substitute $x'=\frac xT$ and pass to Fourier transforms. Then
\begin{equation*}
 B_T=T^{\alpha}\frac 1{2\pi}\int_\R\left(1-e^{-\frac 2{T^\alpha}\abs{x}^\alpha}\right)\hat\varphi(x)\overline{\hat\psi(x)}dx,
\end{equation*}
and it is clear that 
\begin{equation}
 \lim_{T\to \infty}B_T=\frac 1\pi\int_\R\hat\varphi(x)\overline{\hat\psi(x)}\abs x^\alpha dx.\label{e:3.18}
\end{equation}
As $\varphi\in{\cal S}$, we have
\begin{equation}
\left|\varphi(\frac{x+h(x)}{T}+y) -\varphi(\frac{x}{T}+y)\right|\le \frac{C}{T} \label{e:3.19}
\end{equation}
by the Lipschitz property and \eqref{e:3.4}.

As before, denote $\phi(x)=\frac{1}{1+x^2}$. Using
\begin{equation*}
|\varphi(u+v)|\le C(1+u^2)\phi(v),\quad u,v\in\R,
\end{equation*}
we get
\begin{equation}
\label{e:3.20}
|\varphi(\frac{x+h(x)}{T}+y)|\le 2C\phi(\frac xT+y)
\end{equation}
for $T>1$, again by \eqref{e:3.4}.

\eqref{e:3.19} and \eqref{e:3.20} imply, after the usual substitution,
\begin{equation*}
|U_T|\le \frac{C}{T^{1-\alpha}}\int_{\R^2}(\phi(x+y)+|\psi(x+y)|)p_{T^{-\alpha}}(y)dydx=\frac{C}{T^{1-\alpha}}\int_\R(\phi(x)+|\psi(x)|)dx,
\end{equation*}
\begin{align*}
|V_T|&\le \frac{C}{T^{1-\alpha}}\int_{\R^3}(\phi(x+y)+|\psi(x+z)|)p_{T^{-\alpha}}(y)p_{T^{-\alpha}}(z)dydzdx\\
&=\frac{C}{T^{1-\alpha}}\int_\R(\phi(x)+|\psi(x)|)dx.
\end{align*}
Hence
\begin{equation}
\label{e:3.21}
\lim_{T\to\infty}U_T=\lim_{T\to\infty}V_T=0,
\end{equation}
and by \eqref{e:3.14} and \eqref{e:3.18} we obtain
\begin{equation}
\label{e:3.22}
\lim_{T\to\infty}E\langle{Z}_T,\varphi\rangle\langle{Z}_T,\psi\rangle=\frac{1}{\pi}\int_\R\hat{\varphi}(x)\overline{\hat{\psi}(x)}|x|^\alpha dx.
\end{equation}
Now, fix $m$ such that $(1/2)(1-\alpha)m>1$. Similarly as in the proof of Theorem \ref{thm:2.1} (see \eqref{e:3.10}, \eqref{e:3.10a}, \eqref{e:3.10b}), to finish the proof it suffices to show that for any $\varphi\in{\cal S}$,
\begin{equation}
\label{e:3.23}
\lim_{T\to\infty}D_T(\varphi)=0,
\end{equation}
where
\begin{align*}
D_T(\varphi)&=\frac{1}{T^{\frac{(1-\alpha)m}{2}}}\sum_{j\in\Z}E\Bigl|\varphi(\frac{\rho^j}{T}+\frac{\vartheta^j_1}{T})\Bigr|^m\\
&=\frac{1}{T^{\frac{(1-\alpha)m}{2}}}\int_{\R^2}E\Bigl|\varphi(\frac{x+h(x)}{T}+y)\Bigr|^mp_{T^{-\alpha}}(y)dydx
\end{align*}
(cf.\ \eqref{e:3.14}). \eqref{e:3.20} implies that
\begin{equation*}
D_T(\varphi)\le \frac{C}{T^{\frac{(1-\alpha)m}{2}-1}}\int_\R\phi^m(x)dx\longrightarrow 0\quad {\rm as}\ T\to\infty,
\end{equation*}
by our assumption on $m$.

(b) The argument of part (a) can be repeated. The only difference is that the functions of the forms \eqref{e:1.6} - \eqref{e:1.8} are not Lipschitz, so \eqref{e:3.19} does not hold for them. Nevertheless, instead of \eqref{e:3.19}, \eqref{e:3.20} we can use the estimates 
\begin{equation*}
\Bigl|\UNO_{[0,t]}(\frac{x+h(x)}{T}+y)-\UNO_{[0,t]}(\frac xT+y)\Bigr|\le \UNO_{[-\frac 1T,\frac 1T]}(\frac xT+y)+\UNO_{[t-\frac 1T,t+\frac 1T]}(\frac xT+y),
\end{equation*}
\begin{equation*}
\UNO_{[0,t]}(\frac{x+h(x)}{T}+y)\le \UNO_{[-\frac1T,t+\frac 1T]}(\frac xT+y),
\end{equation*}
as well as similar inqualities for $\UNO_{[-t,0]}$. This permits to derive analogues of \eqref{e:3.18}, \eqref{e:3.21} and \eqref{e:3.23}. We omit details.
\qed

\subsection{Proof of Propositions \ref{prop:2.11} and \ref{prop:2.14}}
\label{sec:propos2.11}
Proposition \ref{prop:2.11} can be obtained by the same argument as that applied to derive Theorem 2.1(a) in \cite{functlim2}, therefore we omit the proof. We will sketch the proof of Proposition \ref{prop:2.14}.

Fix $\alpha>1$ and consider $\psi_t=\ind_{[0,t]}-\ind_{[-t,0]}$. For arbitrary $a_1,\ldots,a_n\in \R$, $t_1,\ldots, t_n\ge 0$ denote
\begin{equation}
 \psi=\sum_{k=1}^n a_k\psi_{t_k}.
 \label{e:3.23a}
\end{equation}
It suffices to show that 
\begin{equation}
 \label{e:3.24}
\lim_{T\to \infty} Ee^{i\<Y_T(1),\psi\>}=e^{-\frac 12 \sum_{k,j=1}^n a_ka_jK(t_k,t_j)},
\end{equation}
where $K$ is, up to a constant, the covariance function of nsfBm (see \eqref{e:1.2}). As $\psi$ is an odd function, by \eqref{e:1.9} and the Poisson initial condition we have 
\begin{equation}
 Ee^{i\<Y_T(1),\psi\>}=\exp\left\{\int_\R \left( E e^{i\frac 1{\sqrt T}\int_0^T \psi(x+\vartheta_r)dr}-1 \right)dx\right\},
\label{e:3.25}
\end{equation}
where $\vartheta$ is the standard $\alpha$-stable process.

We start with calculating the covariance
\begin{align}
 E&\int_\R\frac 1T\int_0^T \psi_t(x+\vartheta_r)dr\int_0^T\psi_s(x+\vartheta_u)dudx\notag\\
 &=\frac 2T\int_0^T\int_r^{T}\int_\R\psi_t(x)\T_{u-r}\psi_s(x)dxdudr\notag\\
 &=\frac 1{T\pi}\int_0^T\int_0^{T-r}\int_\R\hat\psi_t(x)e^{-u\abs{x}^\alpha}\overline{\hat\psi_s(x)}dxdudr\notag\\
 &=\frac 4\pi\int_0^1\int_\R\frac{(1-\cos tx)(1-\cos sx)}{\abs{x}^{2+\alpha}}\left(1-e^{-Tr\abs{x}^\alpha}\right)dxdr\notag\\
&\underset{T\to\infty}{\longrightarrow} \frac 4\pi\int_0^1\int_\R\frac{(1-\cos tx)(1-\cos sx)}{\abs{x}^{2+\alpha}}dxdr
=K(t,s)\label{e:3.26}
\end{align}
(cf. \eqref{e:2.17}). Since $\psi$ is odd and $\vartheta$ is symmetric, \eqref{e:3.24} will follow from \eqref{e:3.26} and \eqref{e:3.25} if we show that for any $t\ge 0$
\begin{equation}
 \label{e:3.27}
 \lim_{T\to\infty}R_T=0,
\end{equation}
where 
\begin{equation*}
 R_T=\frac 1{T^2}\int_\R E\left(\int_0^T\psi_t(x+\vartheta_r)dr\right)^4dx.
\end{equation*}
Similarly as before,
\begin{align*}
 \abs{R_T}\le& \frac C{T^2}\int_0^T\int_{r_1}^T\int_{r_2}^T\int_{r_3}^T\int_{\R^3}\abs{\hat\psi_t(x)} e^{-(r_2-r_1)\abs{x}^\alpha} \abs{\hat\psi_t(x-y)} e^{-(r_3-r_2)\abs{y}^\alpha}\\
 &\hskip 2cm\abs{\hat\psi_t(y-z)}e^{-(r_4-r_3)\abs{x}^\alpha}
 \abs{\hat\psi_t(z)}dxdydzdr_4dr_3dr_2dr_1\\
 \le & \frac{C_1}{T}\int_{\R^3}\frac{1-\cos tx}{\abs x} \frac{1-\cos t(x-y)}{\abs {x-y}} \frac{1-\cos t(y-z)}{\abs {y-z}}
  \frac{1-\cos tz}{\abs {z}} \frac{1-e^{-T\abs{y}^\alpha}}{\abs{y}^\alpha}\\
  &\hskip 7cm
  \frac{1}{\abs{x}^\alpha\abs{y}^\alpha}dxdydz.
\end{align*}
Fix any $0<\varepsilon<\frac 1\alpha$. Using $(1-e^{-T\abs{y}^\alpha})\abs{y}^{-\alpha}\le T^{1-\varepsilon}\abs{y}^{-\alpha \varepsilon}$ it is easy to see that
\begin{equation*}
 \abs{R_T}\le \frac {C_2}{T^\varepsilon}\int_{\R^2}\frac{1-\cos tx}{\abs {x}^{1+\alpha}}\frac{1-\cos tz}{\abs {z}^{1+\alpha}}dxdz,
\end{equation*}
hence \eqref{e:3.27} follows. \qed

\subsection*{Proof of Theorem \ref{thm:2.15}}
The idea of the proof is similar to the one used in the proof of Proposition \ref{prop:2.14}. Again, we take $\psi$ of the form \eqref{e:3.23a} and want to prove \eqref{e:3.24}. Let $N^{x,T}$ denote the empirical process of the system started from one particle at point $x$. 

By the Poisson initial condition
\begin{equation}
 \label{e:3.28}
 Ee^{i\<Y_T(1),\psi\>}=\exp\left\{\int_\R \left( E e^{i\frac 1{F_T}\int_0^T \<N^{x,T}_r,\psi\>dr}-1 \right)H_Tdx.
 \right\}
\end{equation}
We will need the following facts about $N^{x,T}$
\begin{align}
 E&\<N^{x,T}_r,\varphi\>=\T_r\varphi(x)
 \label{e:3.29}\\
 E&\<N^{x,T}_r,\varphi_1\>\<N^{x,T}_u,\varphi_2\>
 =\T_r(\varphi_1\T_{u-r}\varphi_2)(x)\notag\\
 &\hskip 3cm+V_T\int_0^r \T_v\left((\T_{r-v}\varphi_1)(\T_{u-v}\varphi_2)\right)(x)dv,\quad r\le u
 \label{e:3.30}\\
 E&\<N^{x,T}_r,\varphi\>^3=\T_r(\varphi^3)+3V_T\int_0^r\T_{r-u}\left(g_{\varphi,T}(u,\cdot)\T_u\varphi\right)(x)dx,
 \label{e:3.31}
\end{align}
where
\begin{equation*}
 g_{\varphi,T}(u,x)=E\<N^{x,T}_r,\varphi\>^2
\end{equation*}
(see \cite{klenke} and \cite{wagneriano} (A.4.3) and (A.4.4)).

\eqref{e:3.29} and the fact that $\psi$ is odd imply that 
\begin{equation*}
  E\int_\R\frac 1{F_T}\int_0^T\<N^{x,T}_r,\psi\>dr H_Tdx=0.
\end{equation*}
\eqref{e:3.24} will follow from \eqref{e:3.28} if we prove that 
\begin{equation}
 A_T(s,t):=\frac{H_T}{F_T^2}E\int_\R\int_0^T\int_0^T\<N^{x,T}_r,\psi_s\>\<N^{x,T}_u,\psi_t\>drdudx
 \underset{T\to\infty}{\longrightarrow} K(s,t)
\label{e:3.32}
\end{equation}
and
\begin{equation}
 B_T(t):=\frac{H_T}{F_T^3}E\int_\R\abs{\int_0^T\<N^{x,T}_r,\psi_t\>dr}^3dx
\underset{T\to\infty}{\longrightarrow} 0
\label{e:3.33}
\end{equation}

By \eqref{e:3.30} and the form of $F_T$, 
\begin{align}
 A_T(s,t)=&\frac 1{TV_T}\int_0^T\int_0^T\int_\R\Bigg( \psi_s(x)\T_{\abs{u-r}}\psi_t)(x)\notag\\
&\left.+ V_T\int_0^{(u\wedge r)} \psi_s(x)(\T_{u+r-2v}\psi_t)(x)dv \right)dxdrdu\notag\\
 =&I_1(T)+I_2(T), \label{e:3.34}
\end{align}
where
\begin{align}
 I_1(T)=&\frac 1{2\pi TV_T}\int_0^T\int_0^T\int_\R\hat\psi_s(x)\overline{\hat \psi_t(x)}e^{-\abs{u-r}\abs{x}^\alpha}dxdrdu\notag \\
 =& \frac 4{\pi V_T}\int_\R \frac{(1-\cos tx)(1-\cos sx)}{x^{2}}\int_0^1\frac{1-e^{-Tr\abs{x}^\alpha}}{\abs{x}^\alpha}drdx\notag\\
 \le & \frac{C(s,t)}{V_T}\to 0, \quad \textrm{as }\ T\to \infty,
 \label{e:3.35}
\end{align}
and
\begin{align}
 I_2(T)&=\frac{1}{2\pi T}\int_\R \hat\psi_s(x)\overline{\hat \psi_t(x)}\int_0^T\int_v^T\int_v^Te^{-(u+v-2r)\abs{x}^\alpha}drdudvdx\notag \\
&=\frac 2{\pi}\int_\R \frac{(1-\cos tx)(1-\cos sx)}{x^{2}}\int_0^1\frac{(1-e^{-Tv\abs{x}^\alpha})^2}{\abs{x}^{2\alpha}}dvdx\notag\\
&\longrightarrow \frac 2\pi \int_\R \frac{(1-\cos tx)(1-\cos sx)}{\abs{x}^{2+2\alpha}}dx, \label{e:3.36}
\end{align}
since $\alpha<\frac 32$.

\eqref{e:3.32} now follows from \eqref{e:3.34}-\eqref{e:3.36} and \eqref{e:2.17}. 

For brevity, denote $f=\abs{\psi_t}$. By \eqref{e:3.31} we have
\begin{align}
 B_T(t)\le & \frac{H_T T^2}{F_T^3}\int_\R \int_0^T E \<N^{x,T}_r,f\>^3drdx\notag\\
 =& II_1(T)+II_2(T)+II_3(T),\label{e:3.37}
\end{align}
where
\begin{align}
 II_1(T)=&\frac{T^{\frac 12}}{H_T^{\frac 12}V_T^{\frac 32}}\int_\R \int_0^T \T_rf^3(x)drdx
 =\frac{T^{\frac 32}}{H_T^{\frac 12}V_T^{\frac 32}}\int_\R f^3(x)dx\to 0
 \label{e:3.38}
\intertext{
as $T\to \infty$ (see \eqref{e:2.16}),}
 II_2(T)=& 3\sqrt{\frac{T}{H_TV_T}}\int_\R\int_0^T\int_0^r \T_uf^2(x)\T_uf(x)dudrdx\notag\\
 \le & C(f)\sqrt{\frac{T}{H_TV_T}} T^2\to 0,\label{e:3.39}
\intertext{as $T\to \infty$,}
 II_3(T)=& \frac{T^{\frac 12}}{H_T^{\frac 12}V_T^{\frac 32}}3V_T^2 \int_\R\int_0^T\int_0^r (\T_uf(x))\int_0^u\T_v(\T_{u-v}f)^2(x)dvdudrdx\notag\\
 \le & C_1(f)\frac{T^{\frac 72}V_T^{\frac 12}}{H_T^{\frac 12}}\to 0
 \label{e:3.40}
\end{align}
as $T\to \infty$ by \eqref{e:2.16}. Thus \eqref{e:3.33} is a consequence of \eqref{e:3.37}-\eqref{e:3.40}. This completes the proof.\qed

%
%

\end{document}